\documentstyle{article}

\begin{document}
\begin{center}
{\Large SOLAR STRUCTURE IN TERMS OF GAUSS'\\
HYPERGEOMETRIC FUNCTION}\\
{\bf H.J. Haubold}\\
UN Outer Space Office, Vienna International Centre, Vienna,
Austria\\
and\\
{\bf A.M. Mathai}\\
Department of Mathematics and Statistics, McGill University,\\
Montreal, Canada
\end{center}

Abstract. Hydrostatic equilibrium and energy conservation determine
the conditions in the gravitationally stabilized solar fusion
reactor. We assume a matter density distribution varying
non-linearly through the central region of the Sun. The analytic
solutions of the differential equations of mass conservation,
hydrostatic equilibrium, and energy conservation, together with the
equation of state of the perfect gas and a nuclear energy
generation rate $\epsilon=\epsilon_0 \rho^nT^m$, are given in terms
of Gauss' hypergeometric function. This model for the structure of
the Sun gives the run of density, mass, pressure, temperature, and
nuclear energy generation through the central region of the Sun.
Because of the assumption of a matter density distribution, the
conditions of hydrostatic equilibrium and energy conservation are
separated from the mode of energy transport in the Sun.
\clearpage
\section{Hydrostatic Equilibrium}
In the following we are concerned with the hydrostatic equilibrium
of the purely
gaseous spherical central region of the Sun generating energy by
nuclear reactions at a certain rate (Chandrasekhar, 1939/1957; Stein, 1966). For this gaseous sphere we
assume that the matter density varies non-linearly from the center
outward, depending on two parameters $\delta$ and $\gamma$,
\begin{equation}
\rho(x)=\rho_cf_D(x),
\end{equation}
\begin{equation}
f_D(x)=[1-x^\delta]^\gamma,
\end{equation}
where $x$ denotes the dimensionless distance variable, $x=r/R_\odot,
0\leq x\leq 1,$ $R_\odot$ is the solar radius, $\delta >0, \gamma>0$
and $\gamma$ is kept a positive integer in the following considerations. The
choice of the density distribution in (1) and (2) reveals
immediately that $\rho(x=0)=\rho_c$ is the central density of the
configuration and $\rho(x=1)=0$ is a boundary condition for hydrostatic equilibrium of the gaseous configuration. For the range
$0\leq x\leq0.3$ the density distribution in (1) and (2)
can be fit numerically to computed data for solar models by chosing
$\delta=1.28$ and $\gamma=10$ (Haubold and Mathai 1994). For these
values of $\delta$ and $\gamma$ the function $f_D(x)$ in (2) is shown in Figure 1. The choice of restricting $x$ to $x\leq 0.3$ is justified by looking at a Standard Solar Model (Bahcall and Pinsonneault, 1992) which shows that $x\leq 0.3$ comprises what is considered to be the gravitationally stabilized solar fusion reactor. More precisely, 95\% of the solar luminosity is produced within the region $x<0.2 (M<0.3M_\odot)$. The half-peak value for the matter density occurs at $x=0.1$ and the half-peak value for the temperature occurs at $x=0.25$. The region $x\leq0.3$ is also the place where the solar neutrino fluxes are generated.
As we are concerned with a spherically symmetrical distribution of
matter, the mass $M(x)$ within the radius $x$ having the density
distribution given in (1) and (2) is
\begin{equation}
M(x)=M_\odot f_M(x),
\end{equation}
\begin{equation}
f_M(x)=\frac{(\frac{3}{\delta}+1)(\frac{3}{\delta}+2)\cdots
(\frac{3}{\delta}+\gamma)}{\gamma!}x^3\, _2F_1(-\gamma, \frac{3}{\delta};\frac{3}{\delta}+1;
x^\delta),
\end{equation}
where $M_\odot$ denotes the solar mass and $_2F_1(.)$ is
Gauss' hypergeometric function (see, for example, Mathai, 1993). 
Equations (3) and (4) are satisfying
the boundary condition $M(x=0)=0$ and determine the central value
$\rho_c$ of the matter density through the boundary condition
$M(x=1)=M_\odot$, where $\rho_c$ depends then only on $\delta$ and
$\gamma$ of the chosen density distribution in (1) and
(2). The function $f_M(x)$ in (4) is shown in Figure
2, using Mathematica (Wolfram, 1991).\par
For hydrostatic equilibrium of the gaseous configuration the
internal pressure needs to balance the gravitational attraction.
The pressure distribution follows by integration of the
respective differential equation for hydrostatic equilibrium,
making use of the density distribution in (1) and the mass
distribution in (3), that is
\begin{equation}
P(x)=\frac{9}{4\pi}G\frac{M_\odot^2}{R_\odot^4}f_P(x),
\end{equation}
\clearpage
\begin{eqnarray}
f_P(x)&=& \left[\frac{(\frac{3}{\delta}+1)(\frac{3}{\delta}+2)\ldots(\frac{3}
{\delta}+\gamma)}{\gamma!}\right]^2\frac{1}{\delta^2}
\sum^\gamma_{m=0}\frac{(-\gamma)_m}{m!
(\frac{3}{\delta}+m)(\frac{2}
{\delta}+m)}\nonumber \\
& & \times
\left[\frac{\gamma!}{(\frac{2}{\delta}+m+1)_\gamma}-x^{\delta
m+2}\,_2F_1(-\gamma,
\frac{2}{\delta}+m;\frac{2}{\delta}+m+1;x^\delta)\right],
\end{eqnarray}
where $G$ is Newton's constant and $_2F_1(.)$ denotes again Gauss'
hypergeometric function (Mathai, 1993).\par
The Pochhammer symbol
$(\frac{2}{\delta}+m+1)_\gamma=\Gamma(\frac{2}{\delta}+m+1+
\gamma)/\Gamma(\frac{2}{\delta}+m+1)$ often appears in series
expansions
for hypergeometric functions. Equations (5) and (6) give the value
of the pressure $P_c$ at the centre of the gaseous configuration
and satisfy the condition $P(x=1)=0.$ The graph of the function
$f_P(x)$ in (6) is shown in Figure 3,using Mathematica (Wolfram, 1991).
\section{Equation of State}
It should be noted that $P(x)$ in (5) denotes the total
pressure of the gaseous configuration, that is the sum of the gas
kinetic pressure and the radiation pressure (according to
Stefan-Boltzmann's law)(Chandrasekhar, 1939/1957; Stein, 1966).
However, the radiation pressure, although the ratio of radiation 
pressure to gas pressure increases
towards the center of the Sun, remains negligibly small in
comparison to the gas kinetic pressure. Thus, Equation (5) can be
considered to represent the run of the gas pressure through the
configuration under consideration. Further, the matter density is
so low that at the temperatures involved the material follows the
equation of state of the perfect gas. Therefore, the temperature
distribution throughout the gaseous configuration is given by
\begin{equation}
T(x)=3\frac{\mu}{kN_A}G\frac{M_\odot}{R_\odot}f_T(x),
\end{equation}
\begin{eqnarray}
f_T(x)&=&\left[\frac{(\frac{3}{\delta}+1)(\frac{3}{\delta}+2)\cdots (\frac{3}{\delta}+\gamma)}{\gamma!}\right]\frac{1}{\delta^2}\frac{1}{[1-x^
\delta]^\gamma}\sum^\gamma_{m=0}
\frac{(-\gamma)_m}{m!(\frac{3}
{\delta}+m)(\frac{2}{\delta}+m)}\nonumber \\
& &
\times\left[\frac{\gamma!}{(\frac{2}{\delta}+m+1)_\gamma}-x^{\delta
m+2}\, _2F_1(-\gamma,
\frac{2}{\delta}+m;\frac{2}{\delta}+m+1;x^\delta)\right],
\end{eqnarray}
where $k$ is the Boltzmann constant, $N_A$ Avogadro's number, $\mu$
the mean molecular weight, and $_2F_1(.)$ Gauss' hypergeometric
function (Mathai, 1993). Equations (7) and (8) reveal the central temperature for
$T(x=0)=T_c$ and satisfy the boundary condition $T(x=1)=0.$ Since
the gas in the central region of the Sun can be treated as
completely ionised, the mean molecular weight $\mu$ is given by
$\mu=(2X+\frac{3}{4}Y+\frac{1}{2}Z)^{-1},$ where $X, Y, Z$ are
relative abundances by mass of hydrogen, helium, and heavy
elements, respectively, and $X+Y+Z=1.$ The run of the function
$f_T(x)$ in (8) is shown in Figure 4, using Mathematica (Wolfram, 1991).
\section{Nuclear Energy Generation Rate}
Hydrostatic equilibrium and energy conservation are determining the
physical conditions in the central part of the Sun. In the
preceeding Sections, the run of density, mass, pressure, and
temperature have been given for a gaseous configuration in
hydrostatic equilibrium based on the equation of state of the
perfect gas (Equations (1)-(8)). In the following, a representation for
the nuclear energy generation rate,
\begin{equation}
\epsilon(\rho, T)=\epsilon_0 \rho^nT^m,
\end{equation}
will be sought which takes into account the above given
distributions of density and temperature and which can be used to
integrate the differential equation of energy conservation 
throughout the gaseous configuration considered in Sections 1 and 2. In
Equation (9), $n$ denotes the density exponent, $m$ the temperature
exponent, and $\epsilon_0$ a positive constant determined by the
specific reactions for the generation of nuclear energy. Using the
equation of state of the perfect gas, Equation (9) can be rewritten
more conveniently,
\begin{eqnarray}
\epsilon(x) & = & \epsilon_0\left(\frac{\mu}{k N_A}\right)^m
\left[\rho(x)\right]^{n-m}\left[P(x)\right]^m \nonumber \\
& & \epsilon_0\left(\frac{\mu}
{kN_A}\right)^m\rho_c^{n-m}P_c^m\left[\frac{P(x)}{P_c}\right]^m
\left[1-x^\delta\right]^{\gamma(n-m)},
\end{eqnarray}
where we note that $0\leq P(x)/P_c\leq 1,$ and $P(x)$ is given in
(5). Subsequently, we can write $P_c$ as follows:
\begin{equation}
P_c=\frac{9}{4\pi}G\frac{M^2_\odot}{R^4_\odot}\left[\frac{(\frac{3}{\delta}+1)(\frac{3}{\delta}+2)\cdots(\frac{3}{\delta}+\gamma)}{\gamma!}\right]^2
\frac{1}{\delta^2}\eta
(\gamma),
\end{equation}
\begin{equation}
\eta(\gamma)=\sum^\gamma_{\nu=0}\frac{(-\gamma)_\nu}{\nu!}
\frac{1}{(\frac{2}{\delta}+\nu)
(\frac{3}{\delta}+\nu)}\frac{\gamma!}{(\frac{2}{\delta}+\nu+1)
\cdots(\frac{2}{\delta}+\nu+\gamma)}.
\end{equation}
Taking the ratio of the pressure $P(x)$ at the location $x$ to the
central value of it one can write
\begin{equation}
\frac{P(x)}{P_c}=1-\frac{1}{\eta(\gamma)}x^2h(x),
\end{equation} 
\begin{eqnarray}
h(x)&=&\sum^\gamma_{m_1=0}\sum^\gamma_{m_2=0}\frac{(-\gamma)_
{m_1}}{m_1!}\frac{(-\gamma)_{m_2}}{m_2!}
\nonumber \\
& & \times
\frac{1}{(\frac{3}{\delta}+m_1)(\frac{2}
{\delta}+m_1+m_2)}x^{\delta(m_1+m_2)}.
\end{eqnarray}
Note that $h(x)$ is a polynomial of degree $2\gamma$ in $x^\delta$.
Denoting the polynomial $h(x)$ by
\begin{equation}
h(x) = a_0+a_1[x^\delta]+a_2[x^\delta]^2+\ldots
+a_{2\gamma}[x^\delta]^{2\gamma},
\end{equation}
we obtain for (13) with a view to (10),
\begin{eqnarray}
\left[\frac{P(x)}{P_c}\right]^m & = &
\left[1-\frac{1}{\eta(\gamma)}x^2h(x)\right]^m \nonumber\\
& = &
\sum^m_{q=0}\frac{(-m)_q}{q!}\left[\frac{1}{\eta(\gamma)}
\right]^qx^{2q}\left[h(x)
\right]^q,
\end{eqnarray}
where we can expand $[h(x)]^q$ by using a multinomial expansion.
That is
\begin{equation}
\left[h(x)\right]^q=\sum^q_{n_0=0}\sum^q_{n_1=0}\ldots\sum^q_
{n_{2\gamma}=0}\frac{q!a_0^{n_0}a_1^{n_1}\ldots
a_{2\gamma}^{n_{2\gamma}}}{n_0!n_1!\ldots
n_{2\gamma}!}\left[x^\delta\right]^{n_1+2n_2+\ldots+(2\gamma)n_
{2\gamma}},
\end{equation}
where $n_0+n_1+\ldots+n_{2\gamma}=q.$ Note also that since $0\leq
x^\delta\leq 1$
we have for $x<1$, in Equation (10),
\begin{equation}
\left[1-x^\delta\right]^{-\gamma(m-n)}=\sum^{\gamma(m-n)}_{s
=0}\frac{[\gamma(m-n)]_s}{s!}
x^{\delta s}.
\end{equation}
For the nuclear energy generation rate in (10) we obtain
finally
\begin{equation}
\epsilon(x)=\epsilon_0\rho_c^nT_c^m f x^{\delta
s+2q+\delta[n_1+2n_2+\ldots +(2\gamma)n_{2\gamma}]},
\end{equation}
where
\begin{eqnarray}
f &=& f(\delta, \gamma, m,n;s,q,n_0,n_1,\ldots,
n_{2\gamma};a_0,a_1,\ldots,a_{2\gamma})\nonumber\\
& &
=\sum^{\gamma(m-n)}_{s=0}\frac{[\gamma(m-n)]_s}{s!}\sum^m_
{q=0}\frac{(-m)_q}{q!}
\left(\frac{1}{\eta(\gamma)}\right)^q\nonumber \\
& & \times \sum^q _{n_0=0}\ldots \sum^q
_{n_{2\gamma}=0}\frac{q!a_0^{n_0}a_1^{n_1}\ldots
a_{2{\gamma}}^{n_{2\gamma}}}{n_0!n_1!\ldots n_{2{\gamma}}!},
\end{eqnarray}
with $n_0+n_1+\ldots +n_{2{\gamma}} = q.$\\
Equation (19) reflects the analytic representation of the nuclear
energy generation rate in (9) taking into account the run
of physical quantities given for the gaseous configuration in
hydrostatic equilibrium considered in Sections 1 and 2.\\
In deriving the explicit dependence of $\epsilon$ on the density
and temperature we have exercised particular care, because the rate
of nuclear energy production is very highly temperature sensitive.
Small changes in the temperature in the central part of the Sun are
adequate to balance large differences in luminosity.
\section{Total Nuclear Energy Generation}
The total net rate of nuclear energy generation is equal to the
luminosity of the Sun, that means the generation of energy by
nuclear reactions in the central part of the Sun has to continually
replenish that energy radiated away at the surface. If $L(x)$
denotes the outflow of energy across the spherical surface at
distance $x$ from the center, then in equilibrium the average
energy production at distance $x$ is
\begin{equation}
L(x)=4\pi R^3_\odot\int_0^x dt t^2\rho (t)\epsilon(t),
\end{equation}
where $\rho(x)$ and $\epsilon(x)$ are given in (1) and
(19), respectively (Chandrasekhar, 1939/1957; Stein, 1966). 
Collecting all factors containing the relative distance
variable $x$, the integral to be evaluated is the following, denoting
it by $g(x)$,
\clearpage
\begin{eqnarray}
g(x) & = & R^3_\odot\int^x_0 dt t^2\rho(t) t^{\delta s+2q+\delta[n_1+2n_2+\ldots+(2\gamma) n_{2\gamma}]}\nonumber\\
&= & \rho_c R^3_\odot \int^x_0 dt t^2[1-t^\delta]^\gamma t^{\delta s+2q+\delta[n_1+2n_2+\ldots+(2\gamma)n_{2\gamma}]}\nonumber\\
&= & \rho_c R^3_\odot \frac{1}{\delta} \int^{x^\delta}_0 dv (1-v)^\gamma v^{s+\frac{1}{\delta}(3+2q)+[n_1+2n_2+\ldots+(2\gamma)n_{2\gamma}]-1}
\end{eqnarray}
by setting $x^\delta = v.$ Equation (22) represents an incomplete beta function which can be written as a series or in terms of a hypergeometric function. Let $s^*=s+\frac{1}{\delta}(3+2q)+n_1+2n_2+\ldots +(2\gamma)n_{2\gamma}$, then we have
\begin{eqnarray}
g(x) & = & \rho_c R^3_\odot \frac{1}{\delta} \int^{x^\delta}_0 dv (1-v)^\gamma v^{s^*-1}\nonumber \\
& = & \rho_c R^3_\odot \frac{1}{\delta}\sum^\gamma _{l=0}\frac{(-\gamma)_l}{l!}\int ^{x^\delta}_0 dv v^{s^*+l-1}\nonumber\\
& = & \rho_c R^3_\odot \frac{1}{\delta}\sum^\gamma _{l=0}\frac{(-\gamma)_l}{l!}\frac{[x^\delta]^{l+s^*}}{l+s^*}\nonumber \\
& = & \rho_c R^3_\odot\frac{1}{\delta}[x^\delta]^{s^*}\sum^\gamma_{l=0}\frac{(-\gamma)_l}{l!}\frac{[x^\delta]^l}{s^*+l}\nonumber\\
& = & \rho_c R^3_\odot\frac{1}{\delta s^*}[x^\delta]^{s^*}\, _2F_1(-\gamma, s^*; s^*+1;x^\delta),
\end{eqnarray} 
where $_2F_1(.)$ is Gauss' hypergeometric function (Mathai, 1993).
Hence
\begin{equation}
L(x)=4\pi \epsilon_0\rho_c^{n+1}T_c^m R^3_\odot\frac{1}{\delta}f\frac{1}{s^*}[x^\delta]^{s^*}\, _2F_1(-\gamma, s^*;s^*+1;x^\delta),
\end{equation}
where f is defined in (20). The luminosity, $L(x=1)=L_\odot$, is given by
\begin{equation}
L_\odot=4\pi \epsilon_0\rho_c^{n+1}T_c^m R^3_\odot\frac{1}{\delta}f\frac{1}{s^*}\frac{\gamma!}{(s^*+1)(s^*+2)\ldots
(s^*+\gamma)},
\end{equation}
using the relation $_2F_1(a,b;c;1)=[\Gamma(c)\Gamma(c-a-b)]/[\Gamma(c-a)\Gamma(c-b)]$
(Mathai, 1993).
\clearpage
\begin{center}
References
\end{center}
Bahcall, J.N. and Pinsonneault, H.M.: 1992, Rev. Mod. Phys. 64,885.\\
Chadrasekhar, S.: 1939/1957, An Introduction to the Study of Stellar Structure,\par 
Dover Publications Inc., New York.\\
Haubold, H.J. and Mathai, A.M.: 1994, in H.J. Haubold and L.I. Onuora (eds.),\par
'Basic Space Science', AIP Conference Proceedings Vol. 320,\par
American Institute of Physics, New York.\\
Mathai, A.M.: 1993, A Handbook of Generalized Special Functions for\par
Statistical and Physical Sciences, Clarendon Press, Oxford.\\
Stein, R.F.: 1966, in R.F. Stein and A.G.W. Cameron (eds.), 'Stellar Evolution',\par
Plenum Press, New York.\\
Wolfram, S.: 1991, Mathematica - A System for Doing Mathematics by\par
Computer, Addison-Wesley Publishing Company Inc., Redwood City.\\
\end{document}